\documentclass[12pt,A4paper,fleqn]{amsart}
\usepackage{mathrsfs}
\usepackage{amsfonts}
\usepackage{txfonts}
\usepackage{amsmath}
\usepackage{amssymb}
\usepackage{amsthm}
\usepackage[pdftex]{graphicx}
\usepackage[toc,page,title,titletoc,header]{appendix}
\usepackage{geometry}
\usepackage{upgreek}
\usepackage[T1]{fontenc}
\usepackage{color}
\usepackage{hyperref}
\usepackage[all]{xy}
\usepackage[french,german,english]{babel}
\usepackage{comment}

\theoremstyle{plain}

\newtheorem*{theo*}{Theorem}

\newtheorem*{lem*}{Lemma}

\theoremstyle{definition}

\theoremstyle{remark}

\newtheorem*{rem*}{Remark}

\newcommand{\R}{\mathbb{R}}

\newcommand{\C}{\mathbb{C}}

\begin{document}

\bibliographystyle{unsrt}

\abstract
In this note, we prove that below the first critical energy level, {a proper combination of the Ligon-Schaaf and Levi-Civita regularization mappings provides} 
a convex {symplectic} embedding of the {energy surfaces} of the planar rotating Kepler problem {into $\R^{4}$ endowed with its standard symplectic structure. A direct consequence is the dynamical convexity of the planar rotating Kepler problem, which has been established in \cite{albers-fish-frauenfelder-koert} by direct computations. This result opens up new approaches to} 
attack the Birkhoff conjecture about the existence of a global surface of section in the restricted planar circular three body problem {using holomorphic curve techniques.}
\endabstract

\title{A convex embedding for the rotating Kepler problem}

\author{Urs Frauenfelder, Otto van Koert, Lei Zhao}


\date\today
\maketitle


\section{Introduction}

The planar rotating Kepler problem is the planar Kepler problem in {a proper uniform-}rotating coordinates. Its Hamiltonian
$$H \colon T^*(\mathbb{R}^2 \setminus \{0\}) \to \mathbb{R}$$
reads
$$H(q,p)=\frac{||p||^2}{2}-\frac{1}{||q||}+(p_1q_2-p_2q_1),$$
{which is the sum of the usual Kepler Hamiltonian in a fixed reference frame with an angular momentum term $(p_1q_2-p_2q_1)$ generating the rotation of the reference frame.}
Note that since the Kepler problem is rotationally invariant{,} the Kepler Hamiltonian
Poisson commutes with the angular momentum. {By c}ompleting the squares{,} we can rewrite the Hamiltonian as
$$
H(q,p)=\frac{1}{2}\Big((p_1+q_2)^2+(p_2-q_1)^2\Big)-\frac{1}{||q||}-\frac{||q||^2}{2}.
$$
The last two terms are referred to as the \emph{effective potential}
$$U(q)=-\frac{1}{||q||}-\frac{||q||^2}{2}$$
so that the Hamiltonian {takes the form}
$$H(q,p)=\frac{1}{2}\Big((p_1+q_2)^2+(p_2-q_1)^2\Big)+U(q).$$
The Hamiltonian of the rotating Kepler problem {is not mechanical, \emph{i.e.} is not a simple sum of}
{kinetic and potential energies}, but contains a twist in the kinetic {part}. 
This twist gives rise to the Coriolis force. {Moreover,}
the effective potential contains {in addition} to the
gravitational potential an elastic term {quadratic in the positions} {which}
gives rise to the centrifugal force. 

Note that the critical values of $H$ correspond to the critical values of $U$. There is precisely one critical value
of $U$ 
$$c_1=-\frac{3}{2}.$$
For an energy value $c$, abbreviate by 
$$\Sigma_c=H^{-1}(c)$$
the energy hypersurface. 
The \emph{Hill's region} of the rotating Kepler problem for the energy value $c$ is the image of the energy hypersurface under the footpoint projection 
$$\pi \colon T^* \mathbb{R}^2 \to \mathbb{R}^2,$$
 i.e.,
$$\mathcal{K}_c=\pi (\Sigma_c)=\Big\{q \in \mathbb{R}^2 \setminus \{0\}: U(q) \leq c\Big\}.$$
For $c<-\tfrac{3}{2}${,} the Hill's region has two components, one bounded and one unbounded
$$\mathcal{K}_c=\mathcal{K}_c^b \cup \mathcal{K}_c^u$$
and we abbreviate by
$$\Sigma_c^b=\Big\{(q,p) \in \Sigma_c: q \in \mathcal{K}_c^b\Big\}$$
the connected component of the energy hypersurface lying over the bounded component of the Hill's region. 
Note that $\Sigma_c^b$ is not compact {due to collisions}. Our main result is
\\ \\
\textbf{Theorem\,A: } \emph{For every $c<-\tfrac{3}{2}$ there exists a 
{2-to-1}-symplectic strictly convex embedding
of $\Sigma_c^b$ to $\mathbb{R}^4$. The image of this embedding smoothly extends to
a closed strictly convex hypersurface in $\mathbb{R}^4$ {over the set of collisions.}}
\\ \\
{At a first glance, one may naturally expect that the embedding obtained by the Levi-Civita regularization suffices. This is unfortunately not the case.}
Indeed, it was shown in \cite[Theorem 1.2]{albers-fish-frauenfelder-koert} that not for all energy levels below the first critical value the Levi-Civita embedding of the
bounded component of the rotating Kepler problem is convex. Instead of that{,} we {obtain our embedding} by first applying the Ligon-Schaaf {regularization mapping}
to the rotating Kepler problem and then apply the Levi-Civita regularization mapping
to the Ligon-Schaaf regularized rotating Kepler problem. 

Ligon and Schaaf discovered their regularization mapping
\cite{ligon-schaaf} in their attempt to understand the symmetries of the Kepler problem by the theory of moment maps. {This regularization mapping can also be thought as a global version of the Delaunay coordinate transformation. }
Nevertheless its mysterious properties still continue to fascinate mathematicians  
{in their efforts to further clarify the situation}, see for example 
 \cite{cushman-duistermaat, heckmann-laat}. 

To prove Theorem\,A{,} we show that the Gauss-Kronecker curvature of the Levi-Civita pull-back of the Ligon-Schaaf regularized rotating Kepler problem is positive. We derive an explicit expression {for this curvature (up to a positive factor)} as a polynomial in three
variables. This polynomial surprisingly factorizes quite well and therefore allows us to estimate the curvature. 

{Since convexity implies dynamical convexity (\cite{HWZ98}, see also Section~\ref{birk}), the first statement of Theorem 1.1 in \cite{albers-fish-frauenfelder-koert}, which states that the bounded component of an energy hypersurface of the rotating Kepler system is dynamically convex, now becomes a direct corollary.}

{Aside from this, w}e do not get {too much} dynamical information {from} Theorem\,A.
Nevertheless Theorem\,A is motivated by dynamical questions. Indeed, the rotating Kepler problem arises as the limit of the restricted three body problem when the mass of one of the primaries goes to zero. In the restricted three body problem{, many} chaotic motion{s occur,} and its dynamics is {still} far from {well-}understood. 
{We hope Theorem\,A will be helpful to understand} a conjecture by Birkhoff about the existence of a global surface of section for the restricted three body problem.
In view of numerical evidence and Theorem\,A it seems reasonable to try to prove the Birkhoff conjecture by verifying convexity of a suitable embedding of the restricted three-body problem into $\C^2$ and applying a result by Hofer, Wysocki and Zehnder.

In Section~\ref{birk} we explain {some old and new progress made towards} the Birkhoff conjecture to further motivate our study. We will also describe some dynamical applications.
In Section~\ref{levi} we prove Theorem\,A.

\section{The Birkhoff conjecture}\label{birk}

Birkhoff writes on page\,328 of his seminal work on the restricted three body problem \cite{birkhoff} 
\begin{quote}
"This state of affairs seems to me to make it probable that the restricted problem of three bodies
admit of reduction to the transformation of a discoid into itself as long as there is a closed oval of zero
velocity about J(upiter) ..."
\end{quote}
In modern mathematical language{, ``}a transformation to a discoid{''} is referred to as the existence of
a disklike global surface of section. The assumption that there is a closed oval of zero velocity 
means that a bounded component of the restricted three body problem for energies below the first critical value
is considered. When Birkhoff referred
to the restricted three body problem, he assumed
that it is regularized {by Levi-Civita regularization.}
\footnote{{Levi-Civita published his paper on the regularization {\cite{levi-civita}} in 1920, but Birkhoff refers to this regularization in his introduction. Goursat already anticipated this in the paper \cite{goursat} published in 1887 .}}
Therefore{,} {when the energy is} below the first critical value{,} the two bounded components are
each diffeomorphic to a three dimensional sphere.

For small energies 
 the Birkhoff conjecture is proved by Conley \cite{conley} and
Kummer \cite{kummer} {for all mass ratios}. 

For sufficiently small mass ratios it was shown by McGehee in \cite{mcgehee} that the Birkhoff conjecture holds true in the connected component around the heavy primary for an arbitrary energy below the first critical value.
That the same results holds as well around the light primary was shown by  Albers, Fish, Frauenfelder, Hofer, and van Koert in 
\cite{albers-fish-frauenfelder-hofer-koert}. 

The proofs by Conley, Kummer and McGehee use{d} perturbative methods{.} 
{I}n contrast to that{,} the proof in \cite{albers-fish-frauenfelder-hofer-koert} is {non-perturbative} in nature. Instead 
it uses global methods of modern symplectic geometry{,}  namely the theory of holomorphic curves in symplectizations{,} due to Hofer \cite{hofer} and Hofer-Wysocki-Zehnder \cite{HWZ98}. Perturbative methods
are only applicable if the system considered is close to a completely integrable system. This holds for small
energy values, where the restricted three body problem is a perturbation of the Kepler problem and for small
mass ratios around the heavy primary where the restricted three body problem is a perturbation of the rotating Kepler problem.  However, for higher energies and higher mass ratios perturbative methods fail. 
Our expectation is that {the use of} holomorphic curves 
{will be} the clue to attack the Birkhoff conjecture. For more information on the relation between the Birkhoff conjecture and the theory of holomorphic curves{,} we refer to
\cite{frauenfelder-koert}. 

In order to construct a
disklike global surface of section via holomorphic curves{,} the question about the existence of a convex embedding
becomes crucial.  
The reasons are as follows. In \cite{HWZ98}{,} Hofer, Wysocki, and Zehnder proved the following result 
\begin{theo*}[Hofer-Wysocki-Zehnder]
Assume that $\Sigma$ is a closed starshaped hypersurface in $\mathbb{R}^4$ endowed with its standard symplectic structure. If $\Sigma$ is dynamically convex, then $\Sigma$ admits a disklike global surface of section. 
\end{theo*}
Here a starshaped hypersurface in $\mathbb{R}^4$ is called dynamically convex if the Conley-Zehnder index of each closed characteristic is at least three. In \cite{albers-frauenfelder-koert-paternain}{,} Albers, Frauenfelder, van Koert and Paternain 
{proved} the following result.
\begin{theo*}[Albers-Frauenfelder-van Koert-Paternain]
The Levi-Civita embedding of each bounded component of the restricted three body problem for energies below the first critical value is a starshaped hypersurface in $\mathbb{R}^4$.
\end{theo*}
In view of 
{these} results{, in order to prove the Birkhoff conjecture, }it suffices to show dynamical convexity for each bounded component of the restricted three body problem for energies below the first critical value. 
{However, t}o check dynamical convexity directly by first determining all closed characteristics{,} and then figuring out their Conley-Zehnder indices is in general not feasible. Instead of that{,} the following result of Hofer-Wysocki-Zehnder from \cite{HWZ98} gives a much more handy approach. 
\begin{theo*}[Hofer-Wysocki-Zehnder]
Assume that $\Sigma \subset \mathbb{R}^4$ is a closed strictly convex hypersurface. Then $\Sigma$ is dynamically convex. 
\end{theo*}
This theorem explains the term "dynamical convex{ity}". While dynamical convexity is a symplectic concept {in the sense that it is preserved under symplectomorphisms,} the notion of convexity is not
. For a convex hypersurface in $\mathbb{R}^4$ there might well be a different starshaped embedding which is not convex. 

 {To the authors' knowledge, t}he question if there exists a dynamically convex, starshaped hypersurface which does not admit a convex embedding at all is {yet open}. In fact {the} dynamical convexity of the rotating Kepler problem was proved already in \cite{albers-fish-frauenfelder-koert} where it was also shown that the Levi-Civita embedding is not {always convex} for all energies. 
 {This} results 
prompted the question if the rotating Kepler problem leads to an example of a dynamically convex hypersurface which does not admit a convex embedding at all. Theorem\,A answers this question in the negative.

We conclude this section with some other 
evidence for the Birkhoff conjecture.
\begin{enumerate}
\item the Levi-Civita regularization of the restricted three-body problem is convex for a wide range of mass ratios and 
 energies.
This 
can be seen rather directly for very small 
energy, nevertheless, 
as was shown in~\cite{albers-fish-frauenfelder-hofer-koert} it also holds for 
energies close to the first critical value provided we are looking at the neighborhood of the primary very light compared to the other primary.

\item numerically, we can verify positivity of the tangential Hessian of the Levi-Civita regularization on a discretization of the energy hypersurface. We then find that the Levi-Civita regularization seems to be convex for all masses $\mu\in (0.01,0.99)$ up to the first critical value.
In fact, convexity of the Levi-Civita regularization only seems to fail very close to $\mu=0$ and $\mu=1$ for energies close to the first critical value.
We emphasize that without precise estimates and careful methods, such as interval arithmetic, such numerical work does not provide a proof, though. 


\item Again on the numerical side, it seems possible to adapt a shooting argument by Birkhoff into a numerical method to find periodic orbits carrying a certain reflection symmetry, to construct a parametrization of the retrograde orbit and a direct orbit. These orbits link like Hopf fibers, and numerical evidence suggests that it is possible to construct a global disklike surface of section from these orbits.



\end{enumerate}

The existence of a global surface of section reveals a lot about the orbit structure; it allows one to study the full dynamics with the globally defined return map, which can be shown to be conjugated to an area-preserving diffeomorphism.
Since a lot is known about the dynamics of such maps, see for instance \cite{franks-handel,lecalvez}, this should provide ample means to better understand
the dynamics.

\section{The Levi-Civita pull-back of the Ligon-Schaaf regularized Kepler Problem}\label{levi}
Set 
$$T=T^{*}\mathbb{S}^{n}=\{(u,v) \in T^{*} \R^{n+1}; \|u\|=1, u\cdot v=0\}$$ and 
$$T^{\times}=\{(u,v) \in T; v \neq 0\}$$ the deleted cotangent bundle of $\mathbb{S}^{n}$, which is sometimes called the Kepler manifold.  Denote by $P_{-}$ the subset of $T^{*} \R^{n}$ with negative Kepler energy $H_{0}(p,q)=\dfrac{\|p\|^{2}}{2}-\dfrac{1}{\|q\|}<0$ and $$T_{-}=\{(u,v) \in T^{\times}; u \neq (0,\cdots,0,1)\}.$$ 

To define the Ligon-Schaaf mapping, we put
\[
\begin{split}
\phi & =-\sqrt{-2K(q,p)}\langle q,p\rangle {,} \\
u &= \left(
\sqrt{-2K(q,p)} \|q\| p,\, \|p\|^2 \|q\|-1
\right) , \\
v &= \left(
-\|q\|^{-1} q+\langle q, p\rangle p,\,\phi
\right)
.
\end{split}
\]
The vectors $u$ and $v$ are orthonormal vectors in $\R^{n+1}$, as can be checked with a direct computation. 
We regard the vector $u$
as the base point in $\mathbb{S}^n$ and the vector $v$ as a unit cotangent vector at $u$.
The Ligon-Schaaf mapping is then given by
\[
\begin{split}
\Phi_{LS}: P_{-} & \longrightarrow T_{-}\\ 
(q,p) & \longmapsto 
\left(
r= (\cos\phi) u+(\sin \phi)v
,s=\frac{1}{\sqrt{-2K(q,p)}}(- (\sin \phi) u + (\cos \phi) v ) 
\,\right),
\end{split}
\]
and it has been shown in \cite{ligon-schaaf,cushman-duistermaat, heckmann-laat} that this map is symplectic with respect to both canonical symplectic structures on cotangent bundles. Furthermore, it transforms $H_{0}(p,q)$ into the ``Delaunay Hamiltonian''
$$H_{k}=-\dfrac{1}{2 \|s\|^{2}}.$$
As we will only study the bounded component of the Hill's region in which all motions are bounded, and thus having negative Keplerian energy, we may well restrict the rotating Kepler problem
$$H=\dfrac{\|p\|^{2}}{2}-\dfrac{1}{\|q\|}+(p_{1} q_{2}-p_{2} q_{1})$$
to $P_{-}$. With the mapping $\Phi_{LS}$, $H$ is transformed into
$$H_{r}=-\dfrac{1}{2 \|s\|^{2}}+(r_{1} s_{2}-r_{2} s_{1})$$
both of which extend smoothly to the north pole $(0,\cdots,0,1)$ of $\mathbb{S}^{n}$ which represents the collisions, and the extensions are thus smoothly defined on $T^{\times}$.

On the other hand, in terms of the semi major axis $a$ and the eccentricity $e$ of the elliptic orbit, the Keplerian energy takes the value $H_{k}=-\dfrac{1}{2 a}$, with the norm of the angular momentum $|p_{1} q_{2}-p_{2} q_{1}|=\sqrt{a} \sqrt{1-e^{2}}$. Moreover, as the bounded component of the Hill's region lies inside the circle $\{\|q\|=1\}$, for all elliptic motions in this component, we have $a <1$. In conclusion, in the bounded component, we have
$$|r_{1} s_{2}- r_{2} s_{1}| \le \|s\| <1.$$

From now on, we shall only consider the planar problem with $n=2$. We have 
$$T^{*} \mathbb{S}^{2}=\{(r_{1}, r_{2}, r_{3}, s_{1}, s_{2}, s_{3}) \in \R^{3} \times \R^{3}; \|r\|=1, r \cdot s=0\}$$
a point $(r_{1}, r_{2}, r_{3}, s_{1}, s_{2}, s_{3})$ which is projected by stereographic projection to
$$(x_{1}, x_{2}, y_{1}, y_{2} )\in \R^{2} \times \R^{2} \cong \C \times \C \ni (x=x_{1} + i x_{2},y=y_{1} + i y_{2})$$
such that
$$s_{1,2}=(\dfrac{\|x\|^{2}+1}{2}) y_{1,2} - Re(\bar{x} y) \cdot x_{1,2}\qquad s_{3}=Re(\bar{x} y)\qquad r_{1,2}=\dfrac{2 x_{1,2}}{\|x\|^{2}+1}, $$
with the north pole projected to ``the point at infinity'' $\infty$.

By calculation, we have
$$\|s\|^{2}=\dfrac{(\|x\|^{2}+1)^{2}}{4} \|y\|^{2} \qquad r_{1} s_{2} -r_{2} s_{1}=x_{1} y_{2}-x_{2} y_{1}.$$

Having in mind the switch in positions and momenta in the Moser regularization, which served as an intermediate step in Heckman-de Laat's interpretation of the Ligon-Schaaf regularization, {we take the following as our Levi-Civita mapping}
$$
L.C.: T^{*} (\C \setminus (0,0)) \to T^{*} \C \qquad (z,w) \mapsto (x=w/ \bar{z}, y=2 z^{2})
$$
which can be extended smoothly to a mapping from $T^{*} \C \setminus \{(0,0)\} \to T^{*} (\C \cup \infty)$.
The pull-back of $H_{r}$ by $L.C.$ thus reads
$$L.C.^{*} H_{r}=-\dfrac{1}{2 (\|w\|^{2}+\|z\|^{2})^{2}}+2 (w_{1} z_{2}-w_{2} z_{1})$$
The corresponding energy level with energy $c$ is therefore
$$\Gamma_{c}=\Bigl\{-\dfrac{1}{2 (\|w\|^{2}+\|z\|^{2})^{2}}+2 (w_{1} z_{2}-w_{2} z_{1})=c\Bigr\}.$$
Note that the angular momentum is $2 (w_{1} z_{2}-w_{2} z_{1})$.
It is envisaged that the set $\Gamma_{c}$ is symmetric with respect to the origin of $\C^{2}$. For $c <- 3/2$, this set is the disjoint union of three components, of which two are unbounded and one is bounded. We would like to understand if the bounded component $\Gamma_{0,c}$ of $\Gamma_{c}$ bounds a convex domain in $\C^{2}$.

In order to show this, we calculate the Gauss-Kronecker curvature of $\Gamma_{0,c}$ and show that this curvature is positive. For this purpose, it is enough to calculate the Hessian of the function 
$$F:=-1+4 (w_{1} z_{2}-w_{2} z_{1})(\|w\|^{2}+\|z\|^{2})^{2}-2 c\, (\|w\|^{2}+\|z\|^{2})^{2}.$$ 
restricted to the tangent space of $\Gamma_{0,c}$ and show that its determinant is positive. The set $\Gamma_{c}$ is just the pre-image $F^{-1}(0)$ of $0$.

To determine the normal direction of points on $\Gamma_{c}$, we calculate the gradient $\nabla F$ of $F$. 
We have 
$$\nabla F=((\|w\|^{2}+\|z\|^{2}) g_{1}, (\|w\|^{2}+\|z\|^{2} )g_{2},  (\|w\|^{2}+\|z\|^{2} ) g_{3},  (\|w\|^{2}+\|z\|^{2} ) g_{4})$$ 
with
\begin{eqnarray*}
\begin{aligned}
& g_{1}= -4 w_{1}^{2} w_{2}+16 w_{1} z_{1} z_{2}-4 w_{2}^{3}-20 w_{2} z_{1}^{2}-4 w_{2} z_{2}^{2}-8 c z_{1}\\
& g_{2}=4 w_{1}^{3}+4 w_{1} w_{2}^{2}+4 w_{1} z_{1}^{2} + 20 w_{1} z_{2}^{2}-16 w_{2} z_{1} z_{2}-8 c z_{2}\\
& g_{3}=20 w_{1}^{2} z_{2}-16 w_{1} w_{2} z_{2}+4 w_{2}^{2} z_{2}+4 z_{1}^{2} z_{2}+4 z_{2}^{3}-8 c w_{1}\\
& g_{4}=-4 w_{1}^{2} z_{1} + 16 w_{1} w_{2} z_{2}-20 w_{2}^{2} z_{1}-4 z_{1}^{3}-4 z_{1} z_{2}^{2}-8 c w_{2}.
\end{aligned}
\end{eqnarray*}
Note that may naturally identify $(g_{1}, g_{2}, g_{3}, g_{4})$ with the quaternion $g:=g_{1}+g_{2} i + g_{3} j + g_{4} k$.  With this identification, we may thus find an orthogonal frame of $T\Gamma_{0,c}$ by (right) multiplications with the quaternions $i,j,k$.
Specifically, we may choose
\begin{eqnarray*}
\begin{aligned}
& v_{1}= (-g_{2}, g_{1}, g_{4}, -g_{3}) \cong g \cdot i 
\\
& v_{2}=(-g_{3}, -g_{4}, g_{1}, g_{2}) \cong g \cdot j\\
& v_{3}=(-g_{4}, g_{3}, -g_{2}, g_{1}) \cong g \cdot k.
\end{aligned}
\end{eqnarray*}
to form a basis of the tangent space at the point $(w_{1}, w_{2}, z_{1}, z_{2})$ provided the gradient 
$$\nabla F=(\Vert w \Vert^2+\Vert z \Vert^2)(g_1,g_2,g_3,g_4)$$
is non-vanishing.
\marginpar{}
We now calculate the determinant $DH$ of the restricted Hessian of $F$ to the tangent spaces of $\Gamma_{0,c}$ and show it is positive. In particular this will show that $\nabla F$ is non-vanishing.
We use the following matrix representation for the tangential Hessian,
$$
DH:= \hbox{Det}\Bigl((v_{1}, v_{2}, v_{3})^{T} \hbox{Hess}(F) (v_{1}, v_{2}, v_{3})\Bigr).
$$
By Maple 18, in expressing the factors in terms of $a=\|w\|^{2}+\|z\|^{2}$, $b=w_{1} z_{2}-z_{1} w_{2}$ and $c$, we find 
$$DH=524288 a^{6} f_{1} f_{2} f_{3} f_{4}^{2}$$
with
\begin{eqnarray*}
\begin{aligned}
& f_{1}=-2 c + a + 4 b\\
& f_{2}=-2 c - a + 4 b\\
& f_{3}=-4 c^{3}+28 b c^{2}-(88 b^{2}-7 a^{2}) c+ 96 b^{3}-15 a^{2 } b\\
& f_{4}= 4 c^{2}- 24 b c + a^{2}+32 b^{2}
\end{aligned}
\end{eqnarray*}

With the conditions $2 | b| \le a <1$ and $-c > 3/2 > (3/2) a$, it is direct to see that 
$$f_{1}>0, f_{2}>0, f_{4}\ge 4(3b-c)^{2}>0.$$

We now show that under the same conditions, we also have $f_{3}>0$. For this, we substitute the relationship $b=\dfrac{1}{4 a^{2}} + c/2$ among $a, b$ and $c$ in the expression of $f_{3}$ and get
$$f_{3}=\dfrac{12 c^{2} a^{4}-2 c a^{8}-15 a^{6}+14 c a^{2}+6}{4 a^{6}}.$$
In which the numerator is a quadratic function in $c$, whose graph is a parabola opening upward with as axis of symmetry the line $c=-\dfrac{7-a^{6}}{12 a^{2}}$. For $a^{2} \ge 7/18$, we have
$$\dfrac{7-a^{6}}{12 a^{2}}<\dfrac{7}{12 a^{2}} < \dfrac{3}{2},$$  and hence this quadratic function is monotonically decreasing for $c<-3/2$. Its evaluation at $c=-3/2$ reads
$$3 a^{8}-15 a^{6}+27 a^{4}-21 a^{2}+6=3(a^{2}-1)^{3} (a^{2}-2)$$
which is clearly positive for $0<a<1$. For $0<a^{2}<7/18$, we find that the evaluation of the denominator of $f_{3}$ at $c=-\dfrac{7-a^{6}}{12 a^{2}}$ reads 
$$-\dfrac{a^{12}}{12}-\dfrac{83 a^{6}}{6}+\dfrac{23}{12},$$
which, as a quadratic equation in $a^{6}$, is seen to be monotonically decreasing when $a^{6}>0$. Moreover its evaluation at $a^{2}=7/18$ is {seen to be positive (approximately $1.1028$)}. This shows that $f_{3}$ is also a positive factor in the factorization of $HD$.

We have thus obtained the conclusion that $\Gamma_{0,c}$ bounds a convex domain for any energy $c$ up to the first critical value $-3/2$. This proves Theorem\,A.

\smallskip
\noindent Urs Frauenfelder, Universit\"at Augsburg: urs.frauenfelder@math.uni-augsburg.de \\
Otto van Koert,   Seoul National University: okoert@snu.ac.kr\\
Lei Zhao,  Chern Institute of Mathematics, Nankai University: l.zhao@nankai.edu.cn


\begin{thebibliography}{99}
 \bibitem{albers-fish-frauenfelder-hofer-koert} P.\,Albers, J.\,Fish, U.\,Frauenfelder, H.\,Hofer,
 O.\,van Koert, \emph{Global surfaces of section in the planar restricted 3-body problem}, Arch.\,Ration.\,Mech.\,Anal.\ \textbf{204}(1) (2012), 273--284.
\bibitem{albers-fish-frauenfelder-koert}P.\,Albers, J.\,Fish, U.\,Frauenfelder, O.\,van Koert, 
\emph{The Conley-Zehnder indices of the rotating Kepler problem}, Math.\,Proc. Cambridge Phil.\,Soc.
\textbf{154}(2)
(2013), 243--260.
 \bibitem{albers-frauenfelder-koert-paternain} P.\,Albers, U.\,Frauenfelder, O.\,van Koert, G.\,Paternain, \emph{Contact geometry of the restricted three-body problem}, Comm.\,Pure Appl.\,Math. \textbf{65}(2) (2012), 229--263.
\bibitem{birkhoff} G.\,Birkhoff, \emph{The restricted problem of
 three bodies}, Rend.\,Circ.\,Matem.\,Palermo \textbf{39} (1915),
 265--334.
\bibitem{conley}C.\,Conley, \emph{On Some New Long Periodic Solutions of the Plane Restricted
 Three Body Problem}, Comm.\,Pure Appl.\,Math. \textbf{16} (1963), 449--467.
\bibitem{cushman-duistermaat}R.\,Cushman, J.J.\,Duistermaat, \emph{A Characterization of the Ligon-Schaaf Regularization Map}, Comm.\,Pure Appl.\,Math. \textbf{50} (1997), 773--787.
\bibitem{franks-handel}
J.\,Franks, M.\,Handel, \emph{Periodic points of Hamiltonian surface diffeomorphisms}
Geom. Topol. 7 (2003), 713--756.
\bibitem{frauenfelder-koert}U.\,Frauenfelder, O.\,van Koert, \emph{The restricted three body problem and holomorphic curves}, book in preparation. 
\bibitem{goursat}E.\,Goursat. \emph{Les transformations isogonales en m\'ecanique,} C. R. Math. Acad. Sci. Paris
, \textbf{108}(1887), 446--450.
\bibitem{heckmann-laat}G.\,Heckman, T.\,de Laat, \emph{On the Regularization of the Kepler Problem},
 Jour.\,Symp.\,Geom. \textbf{10}(3), (2012), 463--473.
\bibitem{hofer}H.\,Hofer, \emph{Pseudoholomorphic curves in symplectisations with application to
the Weinstein conjecture in dimension three}, Invent.\,Math. \textbf{114} (1993), 515--563.
\bibitem{HWZ98} H.\,Hofer, K.\,Wysocki, E.\,Zehnder, {\em The dynamics on a strictly convex energy surface in $\R^4$}, Ann. Math., \textbf{148} (1998) 197--289.
\bibitem{kummer}M.\,Kummer, \emph{On the stability of Hill's Solutions of the Plane Restricted Three Body
Problem}, Amer.\,J.\,Math. \textbf{101}(6), (1979), 1333--1354.
\bibitem{lecalvez} P.\,Le Calvez, \emph{Periodic orbits of Hamiltonian homeomorphisms of surfaces}, Duke Math. J. 133 (2006), no. 1, 125--184.
\bibitem{levi-civita} T.\,Levi-Civita, \emph{Sur la r\'egularisation du probl\`me des trois corps}, Acta. Math. \textbf{42}(1920), 204-219.
\bibitem{ligon-schaaf}T.\,Ligon, M.\,Schaaf, \emph{On the Global Symmetry of the Classical Kepler Problem}, Reports on Math.\,Phys. \textbf{9} (1976), 281--300.
\bibitem{mcgehee}R.\,McGehee, \emph{Some homoclinic orbits for the restricted three-body problem}, Ph.D. Thesis--The University of Wisconsin - Madison. 1969. 63 pp.
\end{thebibliography}
\end{document}